\documentclass[english,a4paper,oneside,12 pt]{article}
\usepackage{amsfonts}
\usepackage{amsmath}
\usepackage[dvips]{graphicx}
\usepackage[english]{babel}
\usepackage{amsmath}
\usepackage{latexsym}
\usepackage{verbatim}
\usepackage{varioref}
\usepackage[latin1]{inputenc}
\usepackage[T1]{fontenc}
\usepackage{array}
\usepackage{delarray}
\usepackage{setspace}
\usepackage{subfigure}
\usepackage{epsfig}

\newenvironment{proof}[1][Proof]{\noindent\\ \textbf{#1} }{\ \rule{0.5em}{0.5em}\\}

\begin{document}

\title{Options on Hedge Funds under the High-Water Mark Rule}
\author{Marc Atlan \\
Laboratoire de Probabilités\\
Université Pierre et Marie Curie \and Hélyette Geman \\
Université Paris Dauphine\\
\& ESSEC Business School \and Marc Yor \\
Laboratoire de Probabilités\\
Université Pierre et Marie Curie}
\date{This draft: February 6 2006}
\maketitle

\begin{abstract}
The rapidly growing hedge fund industry has provided individual and
institutional investors with new investment vehicles and styles of
management. It has also brought forward a new form of performance
contract: hedge fund managers receive incentive fees which are
typically a fraction of the fund net asset value (NAV) above its
starting level - a rule known as \textit{high water mark}.

Options on hedge funds are becoming increasingly popular, in
particular because they allow investors with limited capital to get
exposure to this new asset class. The goal of the paper is to
propose a valuation of plain-vanilla options on hedge funds which
accounts for the high water market rule. Mathematically, this
valuation leads to an interesting use of local times of Brownian
motion. Option prices are numerically computed by inversion of their
Laplace transforms.
\end{abstract}

\noindent \textbf{Keywords: }{Options on hedge funds; High-water
mark; Local time; Excursion theory}

\pagebreak

\section{Introduction}

The term hedge fund is used to characterize a broad class of
"skill-based" asset management firms that do not qualify as mutual
funds regulated by the Investment Company Act of 1940 in the United
States.\ Hedge funds are pooled investment vehicles that are
privately organized, administered by professional investment
managers and not widely available to the general public.\ Due to
their private nature, they carry much fewer restrictions on the use
of leverage, short-selling and derivatives than more regulated
vehicles.

Across the nineties, hedge funds have been embraced by investors
worldwide and are today recognized as an asset class in its own
right. Originally, they were operated by taking a "hedged" position
against a particular event, effectively reducing the overall risk.\
Today, the hedge component has totally disappeared and the term
"hedge fund" refers to any pooled investment vehicle that is not a
conventional fund using essentially long strategies in equity, bonds
and money market instruments.

Over the recent years, multi-strategy funds of funds have in turn
flourished, providing institutional investors with a whole spectrum
of alternative investments exhibiting low correlations with
traditional asset classes.\ In a parallel manner, options on hedge
funds have been growing in numbers and types, offering individual
investors the possibility of acquiring exposure to hedge funds
through a relatively low amount of capital paid upfront at inception
of the strategy.

Hedge funds constitute in fact a very heterogeneous group with
strategies as diverse as convertible arbitrage, global macro or long
short equity.\ In all cases however, common characteristics may be
identified such as long-term commitment of investors, active
management and broad discretion granted to the fund manager over the
investment style and asset classes selected.\ Accordingly, incentive
fees represent a significant percentage of the performance -
typically ranging from 5\% to 20\%.\ This performance is most
generally measured according to the high-water mark rule, i.e.,
using as a reference benchmark the Net Asset Value (NAV) of the fund
at the time of purchase of the shares or options written on the
hedge fund.

So far, the academic literature on hedge funds has focused on such
issues as non-normality of returns, actual realized hedge fund
performance and persistence of that performance.\ Amin and Kat
(2003) show that, as a stand-alone investment, hedge funds do not
offer a superior risk-return profile.\ Geman and Kharoubi (2003)
propose instead the introduction of copulas to better represent the
dependence structure between hedge funds and other asset classes.
Agarwal and Naik (2000) examine whether persistence is sensitive to
the length of the return measurement period and find maximum
persistence at a quarterly horizon.

Another stream of papers has analyzed performance incentives in the
hedge fund industry (see Fung and Hsieh (1999), Brown, Goetzmann and
Ibbotson (1999)).\ However, the high water mark rule specification
has been essentially studied by Goetzman, Ingersoll and Ross (2003).

\bigskip

High-water mark provisions condition the payment of the performance
fee upon the hedge fund Net Asset Value exceeding the entry point of
the investor.\ Goetzmann et al examine the costs and benefits to
investors of this form of managers' compensation and the
consequences of thess option-like characteristics on the values of
fees on one hand, investors' claims on the other hand.\ Our
objective is to pursue this analysis one step further and examine
the valuation of options on hedge funds under the high-water mark
rule.\ We show that this particular compound option-like problem may
be solved in the Black-Scholes (1973) and Merton (1973) setting of
geometric Brownian motion for the hedge fund NAV by the use of Local
times of Brownian motion.

The remainder of the paper is organized as follows: Section II
contains the description of the Net Asset Value dynamics, management
and incentive fees and the NAV option valuation.\ Section II also
extends the problem to a moving high water mark.\ Section III
describes numerical examples obtained by inverse Laplace transforms
and Monte Carlo simulations.\ Section IV contains concluding
comments.

\section{The High-Water Mark Rule and Local Times}

\subsection*{A.\indent Modeling the High-Water Mark}

We work in a continuous-time framework and assume that the fund Net
Asset\ Value (NAV) follows a lognormal diffusion process. This
diffusion process will have a different starting point for each
investor, depending on the time she entered her position.\ This
starting point will define the high water mark used as the benchmark
triggering the performance fees discussed throughout the paper.

We follow Goetzmann, Ingersoll and Ross (2003) in representing the performance fees in the following form%
\begin{equation}
f\left( S_{t}\right) =\mu a~\mathbf{1}_{\{S_{t}>H\}}
{\label{perffeesf}}
\end{equation}
\noindent where $S_{t}$ denotes the Net Asset Value at date t, $\mu
$ is a mean NAV return statistically observed, $a$ is a percentage
generally comprised between 5\% and 20\% and $H=S_{0}$ denotes the
market value of the NAV as observed at inception of the option
contract.

We consider $(\Omega ,\mathcal{F},\{\mathcal{F}_{t},t\geq
0\},\mathbb{P}_0)$ a filtered probability space where
$(B_{t})_{t\geq 0}$ is an $\{\mathcal{F}_{t},t\geq 0\}$ Brownian
motion.

We now consider an equivalent measure $\mathbb{Q}$ under which the
Net Asset Value dynamics $(S_{t})_{t\geq 0}$ satisfy the stochastic
differential equation:
\begin{equation}
\frac{dS_{t}}{S_{t}}=(r+\alpha -c-f(S_{t}))dt+\sigma
dW_{t}{\label{fundSDE}}
\end{equation}

\noindent and the instantaneously compounding interest rate $r$ is
supposed to be constant. $\alpha $ denotes the excess return on the
fund's assets and is
classically defined by%
\begin{equation*}
\alpha =\mu -r-\beta \left( r_{m}-r\right)
\end{equation*}
\noindent where $r_{m}$ is the expected return on the market
portfolio.\ Hence, the "risk-neutral" return on the fund NAV is
equal to $(r+\alpha )$\footnote{\renewcommand{\baselinestretch}{2.0}
Our claim is that the measure $\mathbb{Q}$
 incorporates the price of market risk as a whole but not the excess performance - the fund "alpha" - achieved by the manager
 through the selection of specific securities at a given point in time. This view is in agreement with the footnote 6 in Goetzmann, Ingersoll and Ross (2003)}; $~\sigma ~$%
denotes the NAV volatility.

The management fees paid regardless of the performance are
represented by a constant fraction $c$ (comprised in practice
between 0.5\% and 2\%) of the Net Asset Value. We represent the
incentive fees as a deterministic function $f$ of the current value
$S_{t}$ of the NAV, generally chosen according to the high water
mark rule defined in equation (1). We can note that management fees
have the form of the constant dividend payment of the Merton (1973)
model while performance fees may be interpreted as a more involved
form of dividend paid to the manager.

Because of their central role in what follows, we introduce the
maximum and the minimum processes of the Brownian motion $B$, namely
\begin{equation*}
M_{t}=\underset{s\leq t}{\sup }B_{s},~I_{t}=\underset{s\leq t}{\inf
}B_{s}
\end{equation*}
\noindent as well as its local time at the level $a$, $a\in
~\mathbb{R}$
\begin{equation*}
L_{t}^{a}=\underset{\epsilon \rightarrow 0}{\lim }~\frac{1}{2\epsilon }%
\int_{0}^{t}\mathbf{1}_{\{|B_{s}-a|\leq \epsilon \}}ds
\end{equation*}

We also consider $A_{t}^{(a,+)}=\int_{0}^{t}\mathbf{1}_{\{B_{s}\geq a\}}ds$ and $%
A_{t}^{(a,-)}=\int_{0}^{t}\mathbf{1}_{\{B_{s}\leq a\}}ds$,
respectively denoting the time spent in $[a;\infty \lbrack $ and the
time spent in $]-\infty ;a]$ by the Brownian motion up to time
t.\newline For simplicity, we shall write $L_{t}=L_{t}^{0}$,
$A_{t}^{+}=A_{t}^{(0,+)}$ and $A_{t}^{-}=A_{t}^{(0,-)}$ the
corresponding quantities for $a=0$.

\bigskip

\noindent In order to extend our results to different types of
incentive fees, we do not specify the function $f$ but only assume
that it is a continuous, bounded, increasing and positive function
satisfying the following conditions:
\begin{equation*}
f(0)=0,\indent\lim_{x\rightarrow \infty }f(x)<+\infty
\end{equation*}

\newtheorem{fund}{Proposition}[section]
\begin{fund}
There exists a unique solution to the stochastic differential
equation
\begin{equation*}
\frac{dS_t}{S_t}=(r+\alpha-c-f(S_t))dt+\sigma dW_t
\end{equation*}
\end{fund}

\begin{proof}
Let us denote $Y_{t}=\frac{\ln (S_{t})}{\sigma }$. Applying Itô's
formula, we see that the process $Y_{t}$ satisfies the equation
\begin{equation*}
dY_{t}=dW_{t}+\psi (e^{\sigma Y_{t}})dt
\end{equation*}
\noindent where $\psi (x)=r-\frac{\sigma
^{2}}{2}+\alpha-c-f(x)$.\newline $f$, hence $\psi $ is a Borel
bounded function; consequently, we may apply Zvonkin (1974) theorem
and obtain strong existence and pathwise uniqueness of the solution
of equation ({\ref{fundSDE}}).\newline
We recall that Zvonkin theorem establishes that for every bounded Borel function $%
\xi$, the stochastic differential equation
\begin{equation*}
dZ_t=dW_t+\xi(Z_t)dt
\end{equation*}
has a unique solution which is strong, i.e.: in this case, the
filtration of $Z$ and $W$ are equal.
\end{proof}\\
\noindent Integrating equation ({\ref{fundSDE}}), we observe that
this unique
solution can be written as%
\begin{equation*}
S_{t}=S_{0}\exp \left( \left( r+\alpha-c-\frac{\sigma ²}{2}\right)
t-\int_{0}^{t}f(S_{u})du+\sigma W_{t}\right)
\end{equation*}\\
We now seek to construct a new probability measure $\mathbb{P}$
under which the expression of $S_{t}$ reduces to
\begin{equation}
S_{t}=S_{0}\exp (\sigma \widetilde{W}_{t}){\label{rem}}
\end{equation}
\noindent where $\widetilde{W}_{t}$ is a $\mathbb{P}$ standard
Brownian motion.

\newtheorem{girsanov}[fund]{Proposition}
\begin{girsanov}
There exists an equivalent martingale measure $\mathbb{P}$ under
which the Net Asset Value dynamics satisfy the stochastic
differential equation
\begin{equation}
\frac{dS_t}{S_t}=\frac{\sigma^2}{2} dt +\sigma d\widetilde{W}_t
\end{equation}
where
\begin{equation}
\mathbb{Q}_{|\mathcal{F}_t}=Z_t\cdot\mathbb{P}_{|\mathcal{F}_t}
{\label{proba}}
\end{equation}
\[Z_t=\exp\big(\int_0^t
\big(b-\frac{f(S_u)}{\sigma}\big)d\widetilde{W}_u-\frac{1}{2}\int_0^t
\big(b-\frac{f(S_u)}{\sigma}\big)^2du\big)\] and
\[b=\frac{r+\alpha-c-\frac{\sigma^2}{2}}{\sigma}\]
\end{girsanov}

\begin{proof}
Thanks to Girsanov theorem (see for instance McKean (1969) and Revuz
and Yor (2005)) we find that under the probability measure
$\mathbb{\mathbb{P}}$,\newline
$\widetilde{W}_t=W_t+\int_0^t du \big(b-\frac{f(e^{\sigma Y_u})}{\sigma}%
\big) $ is a Brownian motion, which allows us to conclude.
\end{proof}

\subsection*{B.\indent Building the Pricing Framework}

For practical purposes, the issuer of the call is typically the
hedge fund itself, hence hedging arguments allow to price the option
as the expectation (under the right probability measure) of the
discounted payoff. More generally, a European-style hedge fund
derivative with maturity $T>0$ is defined by its payoff $F:
\mathbb{R}_{+}\longrightarrow \mathbb{R}_{+}$ and the valuation of
the option reduces to computing expectations of the following form:
\begin{equation*}
V_{F}(t,S,T)=e^{-r(T-t)}\mathbb{E^{Q}}\big[F(S_{u};u\leq T)\big|\mathcal{F}%
_{t}]{\label{fundprice}}
\end{equation*}%
For the case where the valuation of the option takes place at a date
$t=0$, we denote $V_{F}(S,T)=V_{F}(0,S,T)$. We can observe that we
are in a situation of complete markets since the only source of
randomness is the Brownian motion driving the NAV dynamics.\newline
\newtheorem{highwat2}[fund]{Proposition}
\begin{highwat2}
For any payoff $F$ that can be written as an increasing function of
the stock price process, the option price associated to the above
payoff is an increasing function of the high-water mark level.
\end{highwat2}

\begin{proof}
This result is quite satisfactory from a financial perspective.
Mathematically, it may be deduced from the following result :\\
Let us consider the solutions $(S^1, S^2)$ of the pair of stochastic
differential equations :
\begin{eqnarray*}
dS^1_{t}&=&b^1(S^1_t)dt+\sigma S^1_t dW_{t}\\
dS^2_{t}&=&b^2(S^2_t) dt+\sigma S^2_t dW_{t}
\end{eqnarray*}
where
\begin{eqnarray*}
b^1(x)&=&(r+\alpha -c-\mu
a~\mathbf{1}_{\{x>H\}})x\\
b^2(x)&=&(r+\alpha -c-\mu a~\mathbf{1}_{\{x>H'\}})x
\end{eqnarray*}
with $H>H'$ and $S^1_0=S^2_0$ a.s.\\
We may apply a comparison theorem since $b^1$ and $b^2$ are bounded
Borel functions and $b^1\geq b^2$ everywhere, obtain that
\begin{equation*}
\mathbb{P}[S^1_t\geq S^2_t; \forall t\geq 0]=1
\end{equation*}
and then conclude.
\end{proof}

If we consider a call option and a put option with strike $K$ and maturity $%
T $, we observe the following call-put parity relation:
\begin{equation}
C_{0}(K,T)-P_{0}(K,T)=\mathbb{E^{Q}}[e^{-rT}S_{T}]-Ke^{-rT}
\end{equation}

\noindent We now wish to express the exponential
$(\mathcal{F}_t,\mathbb{P})$-martingale $Z_{t}$ featured in
(\ref{proba}) in terms of well-known processes in order to be able
to obtain closed-form pricing formulas.

\newtheorem{expgirs}[fund]{Lemma}
\begin{expgirs}
Let us define $d_H$, $\lambda$, $\alpha_{+}$, $\alpha_{-}$ and
$\phi$ as follows:
\[d_H=\frac{\ln(\frac{H}{S_0})}{\sigma},\indent \lambda=\frac{\mu a}{2\sigma}\]
\[\alpha_{+}=2\lambda^2+\frac{b^2}{2}-2\lambda b ,\indent \alpha_{-}=\frac{b^2}{2}\]
\[\phi(x)=e^{b x-2\lambda (x-d_H)_{+}}\]
We then obtain:
\begin{equation}
Z_t=e^{2\lambda(-d_H)_{+}}\phi(\widetilde{W}_t)\exp(\lambda
L^{d_H}_t)\exp(-\alpha_{+}A_t^{(d_H,+)}-\alpha_{-}A_t^{(d_H,-)})
\end{equation}
\end{expgirs}

\begin{proof}
The proof of this proposition is based on the one hand on the Tanaka
formula which, for a Brownian motion $B$ and any real number $a$,
establishes that
\begin{equation*}
(B_{t}-a)_{+}=(-a)_{+}+\int_{0}^{t}dB_{s}\mathbf{1}_{\{B_{s}>a\}}+\frac{1}{2}%
L_{t}^{a}
\end{equation*}%
On the other hand, we can rewrite
\begin{equation*}
f(S_{t})=\mu a\mathbf{1}_{\{\widetilde{W}_{t}>d_{H}\}}
\end{equation*}%
Observing that $A_{t}^{(d_{H},+)}+A_{t}^{(d_{H},-)}=t$ leads to the
result.
\end{proof}\\
\noindent From the above lemma, we obtain that:
\begin{eqnarray*}
V_{F}(S,T) &=&e^{-rT}\mathbb{E^{P}}\big[Z_{T}F(S_{u};u\leq T)\big] \\
&=&e^{-rT+2\lambda (-d_{H})_{+}}\mathbb{E^{P}}\big[\phi (W_{T})\exp
(\lambda L_{T}^{d_{H}}-\alpha _{+}A_{T}^{(d_{H},+)}-\alpha
_{-}A_{T}^{(d_{H},-)})F(S_{0}e^{\sigma \widetilde{W}_{u}};u\leq
T)\big]
\end{eqnarray*}
The price of a NAV call option is closely related to the law of the
triple $(W_{t},L_{t}^{a},A_{t}^{(a,+)})$. Karatzas and Shreve (1991)
have extensively studied this joint density for $a=0$ and obtained
in particular
the following remarkable result\\

\newtheorem{kars}[fund]{Proposition}
\begin{kars}
For any positive $t$ and $b$, $0<\tau<t$, we have
\begin{eqnarray*}
\mathbb{P}[W_t\in dx; L_t\in db,A^{+}_t\in d\tau]&=&f(x,b;t,\tau)~dx  ~db ~d\tau;\indent x>0\\
&=&f(-x,b;t,-\tau)~dx  ~db  ~d\tau;\indent x<0
\end{eqnarray*}
where
\begin{equation*}
f(x,b;t,\tau)=\frac{b(2x+b)}{8\pi\tau^{\frac{3}{2}}(t-\tau)^{\frac{3}{2}}}\exp{\bigg(-\frac{b^2}{8(t-\tau)}-\frac{(2x+b)^2}{8\tau}\bigg)}
\end{equation*}
\end{kars}

This formula could lead to a computation of the option price based
on a multiple integration but it would be numerically intensive;
moreover, obtaining an analytical formula for the triple integral
involved in the option price seems quite unlikely. We observe
instead that in the above density $f$, a convolution product
appears, which leads us to compute either Fourier or Laplace
transforms. We are in fact going to compute the Laplace transform
with respect to time to maturity of the option price. This way to
proceed is mathematically related to the Karatzas and Shreve result
in Proposition II.5. In the same way, we can notice that the Laplace
transform exhibited by Geman and Yor (1996) for the valuation of a
Double Barrier option is related to the distribution of the triple
$(W_t,M_t,I_t)$ Brownian motion, its maximum and minimum used by
Kunitomo and Ikeda (1992) for the same pricing problem. The formulas
involved in the NAV call price rely on the following result which
may be obtained from Brownian excursion theory:

\newtheorem{excur}[fund]{Proposition}
\begin{excur}
Let $W_t$ be a standard Brownian motion, $L_t$ its local time at
zero, $A_t^{+}$ and $A_t^{-}$ the times spent positively and
negatively until time $t$. For any function $h\in L^1(\mathbb{R})$,
the Laplace transform of the quantity
$g(t)=\mathbb{E}\big[h(W_t)\exp(\lambda L_t)\exp(-\mu A_t^{+}-\nu
A_t^{-})\big]$ has the following analytical expression
\begin{equation*}
\int_0^{\infty} dt
e^{-\frac{\theta}{2}t}g(t)=2\frac{\bigg(\int_0^{\infty}dx
e^{-x\sqrt{\theta+2\mu}}h(x)+\int_0^{\infty}dx
e^{-x\sqrt{\theta+2\nu}}h(-x)\bigg)}{\sqrt{\theta+2\mu}+\sqrt{\theta+2\nu}
-2\lambda}
\end{equation*}
for $\theta$ large enough to ensure positivity of the denominator.
\end{excur}

\begin{proof}
See the Appendix for details. The result is rooted in the theory of
excursions of the Brownian motion.
\end{proof}

\subsection*{C.\indent Valuation of the Option at Inception of the Contract}

In this section, we turn to the computation of the price of a
European call option written on a Hedge Fund NAV under the
high-water mark rule. Consequently, the payoff considered is the
following:
\begin{equation}
F(S_{u};u\leq T)=(S_{T}-K)_{+}
\end{equation}%
or, in a more convenient way for our purpose
\begin{equation*}
F(\widetilde{W}_{u};u\leq T)=(S_{0}\exp (\sigma
\widetilde{W}_{T})-K)_{+}
\end{equation*}%
At inception of the contract, the high-water mark that is chosen is
the spot price, hence $H=S_{0}$ and $d_{H}=0$. This specific
framework allows us to use fundamental results on the joint law of the triple $%
(B_{t},L_{t}^{0},A_{t}^{+})$ presented in Proposition II.6. We write
the European call option price as follows
\begin{equation*}
C(0,S_0)=e^{-rT}\mathbb{E^{P}}\big[h(\widetilde{W}_{T})\exp (\lambda
L_{T}-\alpha _{+}A_{T}^{+}-\alpha _{-}A_{T}^{-})\big]
\end{equation*}%
where $h(x)=(S_{0}e^{\sigma x}-K)_{+}e^{bx-2\lambda
(x)_{+}}$.\newline

We now compute the Laplace transform in time to maturity of the
European call option on the NAV of an Hedge Fund, that is to say the
following quantity:
\begin{eqnarray*}
\forall \theta \in \mathbb{R}_{+}\indent I(\theta )&=&\int_{0}^{\infty }dte^{-%
\frac{\theta}{2}t}e^{-rt}\mathbb{E^{Q}}[(S_{t}-K)_{+}]\\
&=&\int_{0}^{\infty }dte^{-(\frac{\theta}{2}+r)t}\mathbb{E^{P}}%
[Z_{t}(S_{t}-K)_{+}]
\end{eqnarray*}%
\newline
\newtheorem{price}[fund]{Lemma}
\begin{price}
The Laplace transform with respect to time to maturity of a call
option price  has the following analytical expression:
\begin{equation}
I(\theta)=2\frac{\bigg(\int_0^{\infty}dx
e^{-x\sqrt{\theta+2(r+\alpha_{+})}}h(x)+\int_0^{\infty}dx
e^{-x\sqrt{\theta+2(r+\alpha_{-})}}h(-x)\bigg)}{\sqrt{\theta+2(r+\alpha_{+})}+\sqrt{\theta+2(r+\alpha_{-})}
-2\lambda}
\end{equation}
where $h(x)=e^{bx-2\lambda x_{+}}(S_0 e^{\sigma x}-K)_{+}$
\end{price}

\begin{proof}
We obtain from Lemma II.4 that:
\begin{equation*}
\mathbb{E^{P}}[Z_{t}(S_{t}-K)_{+}]=\mathbb{E}\big[h(\widetilde{W}_{t})\exp
(\lambda L_{t})\exp (-\alpha _{+}A_{t}^{+}-\alpha
_{-}A_{t}^{-})\big]
\end{equation*}%
where:
\begin{equation*}
h(x)=e^{bx-2\lambda x_{+}}(S_{0}e^{\sigma x}-K)_{+}
\end{equation*}%
Then, using Proposition II.6, we are able to conclude.
\end{proof}

This lemma leads us to compute explicit formulas for the Laplace
transform
of a call option that is in-the-money $(S_{0}\geq K)$ at date 0 and out-of-the-money $%
(S_{0}<K)$ that we present in two consecutive propositions.\\
\newtheorem{otm}[fund]{Proposition}
\begin{otm}
For an out-of-the-money call option ($S_0\leq K$), the Laplace
transform of the price is given by the following formula:
\[I(\theta)=\frac{N(\theta)}{D(\theta)}\]
where
\[\theta>(\sigma+b-2\lambda)^2-2(r+\alpha_{+})\] and
\begin{eqnarray*}
D(\theta)&=&\frac{\sqrt{\theta+2(r+\alpha_{+})}+\sqrt{\theta+2(r+\alpha_{-})}-2\lambda}{2}\\
N(\theta)&=&\frac{S_0}{\sqrt{\theta+2(r+\alpha_{+})}+2\lambda-\sigma-b}\bigg(\frac{S_0}{K}\bigg)^\frac{\sqrt{\theta+2(r+\alpha_{+})}+2\lambda-\sigma-b}{\sigma}\\&&-\frac{K}{\sqrt{\theta+2(r+\alpha_{+})}+2\lambda-b}\bigg(\frac{S_0}{K}\bigg)^\frac{\sqrt{\theta+2(r+\alpha_{+})}+2\lambda-b}{\sigma}
\end{eqnarray*}
\end{otm}

\begin{proof}
Keeping the notation of Proposition II.6, we can write
\begin{equation*}
\forall x>0,\indent h(x)=(S_{0}e^{\sigma x}-K)\mathbf{1}_{\{x\geq \frac{1}{%
\sigma }\ln (\frac{K}{S_{0}})\}}e^{(b-2\lambda )x}\indent and\indent
h(-x)=0
\end{equation*}%
and then by simple integration, obtain the stated formula.
\end{proof}

\newtheorem{itm}[fund]{Proposition}
\begin{itm}
For an in-the-money call option ($S_0\geq K$), the Laplace transform
of the price is given by the following formula:
\begin{equation*}
I(\theta)=\frac{N_1(\theta)+N_2(\theta)}{D(\theta)}
\end{equation*}
where $\theta>(\sigma+b-2\lambda)^2-2(r+\alpha_{+})$ and
\begin{eqnarray*}
D(\theta)&=&\frac{\sqrt{\theta+2(r+\alpha_{+})}+\sqrt{\theta+2(r+\alpha_{-})}-2\lambda}{2}\\
N_1(\theta)&=&\frac{S_0}{\sqrt{\theta+2(r+\alpha_{+})}+2\lambda-\sigma-b}-\frac{K}{\sqrt{\theta+2(r+\alpha_{+})}+2\lambda-b}\\
N_2(\theta)&=&\frac{S_0}{\sqrt{\theta+2(r+\alpha_{-})}+\sigma+b}\bigg(1-\big(\frac{K}{S_0}\big)^{\frac{\sqrt{\theta+2(r+\alpha_{-})}+\sigma+b}{\sigma}}\bigg)\\
&&-\frac{K}{\sqrt{\theta+2(r+\alpha_{-})}+b}\bigg(1-\big(\frac{K}{S_0}\big)^{\frac{\sqrt{\theta+2(r+\alpha_{-})}+b}{\sigma}}\bigg),
\end{eqnarray*}
$\alpha_{-}$ and $\alpha_{+}$ being defined in Lemma II.4 .
\end{itm}

\begin{proof}
We have
\begin{equation*}
\forall x>0,\indent h(x)=(S_0 e^{\sigma x}-K)e^{(b-2\lambda)x}\indent and%
\indent h(-x)=(S_0 e^{-\sigma x}-K)\mathbf{1}_{\{x\leq\frac{1}{\sigma}\ln(%
\frac{S_0}{K})\}}e^{-b x}
\end{equation*}
and as in the previous proposition, the Laplace transform is
derived.
\end{proof}\\
As a side note, we observe that the case $K=0$ provides the Laplace
transform of the $t-$maturity forward contract written on the NAV at
date $0$
\begin{equation*}
\int_0^{\infty} dt e^{-\frac{\theta}{2}t}\mathbb{E^P}[e^{-rt}S_t]=2\frac{%
\frac{S_0}{\sqrt{\theta+2(r+\alpha_{+})}+2\lambda-\sigma-b}+\frac{S_0}{%
\sqrt{\theta+2(r+\alpha_{-})}+\sigma+b}}{\sqrt{\theta+2(r+\alpha_{+})}+%
\sqrt{\theta+2(r+\alpha_{-})}-2\lambda}
\end{equation*}
where $\theta>(\sigma+b-2\lambda)^2-2(r+\alpha_{+})$.\\

It is satisfactory to check that by choosing $a=0$, $\alpha=0$, we
obtain the Laplace transform of a European call option on a
dividend-paying stock with a continuous dividend yield $c$ whose
dynamics satisfy as in Merton (1973), the equation
\begin{equation*}
\frac{dS_t}{S_t}=(r-c)~dt+\sigma ~dW_t
\end{equation*}
This Laplace transform is derived from Proposition II.8 for an
out-of-the-money call option and from Proposition II.9 for an
in-the-money call option.

\subsection*{D.\indent Valuation during the lifetime of the Option}

Evaluating at a time $t$ a call option on a hedge fund written at
date $0$ implies that we are in the situation where $d_{H}=\frac{1}{%
\sigma }\ln (\frac{H}{S_{t}})$ may be different from $0$. Since the
solution of the stochastic differential equation driving the Net
Asset Value is a Markov process, the evaluation of the option at
time $t$ only depends on the value of the process at time $t$ and on
the time to maturity $T-t$. Hence, we need to compute the following
quantity
\begin{equation*}
C(t,S_t)=\mathbb{E^{Q}}\big[e^{-r(T-t)}(S_T-K)_{+}|\mathcal{F}_t\big]
\end{equation*}
Given the relationship between $\mathbb{P}$ and $\mathbb{Q}$, we can
write
\begin{equation*}
C(t,S_t)=e^{-r(T-t)}e^{2\lambda (-d_{H})_{+}}\mathbb{E^{P}}\big[h(%
\widetilde{W}_{T-t})\exp (\lambda L_{T-t}^{d_{H}}-\alpha
_{+}A_{T-t}^{(d_{H},+)}-\alpha _{-}A_{T-t}^{(d_{H},-)})\big]
\end{equation*}%
where $h(x)=e^{bx-2\lambda (x-d_{H})_{+}}(S_{t}e^{\sigma x}-K)_{+}$
\newline
Because of the importance of the level $d_{H}$ in the computations,
we introduce the stopping time $\tau _{d_{H}}=\inf \{t\geq
0;\widetilde{W}_{t}=d_{H}\}$ and split the problem into the
computation of the two following quantities:
\begin{eqnarray*}
C_1&=&e^{-r(T-t)}e^{2\lambda (-d_{H})_{+}}\mathbb{E^{P}}\big[\mathbf{1}_{\{\tau _{d_{H}}>T-t\}}h(\widetilde{W%
}_{T-t})\exp (\lambda L_{T-t}^{d_{H}}-\alpha
_{+}A_{T-t}^{(d_{H},+)}-\alpha
_{-}A_{T-t}^{(d_{H},-)})\big] \\
\textrm{and}\\
C_2&=&e^{-r(T-t)}e^{2\lambda (-d_{H})_{+}}\mathbb{E^{P}}\big[\mathbf{1}_{\{\tau _{d_{H}}<T-t\}}h(\widetilde{W%
}_{T-t})\exp (\lambda L_{T-t}^{d_{H}}-\alpha
_{+}A_{T-t}^{(d_{H},+)}-\alpha _{-}A_{T-t}^{(d_{H},-)})\big]
\end{eqnarray*}%

In order to compute $C_1$, we introduce for simplicity $s=T-t$ and
obtain
\begin{eqnarray*}
\mathbb{E^P}\big[\mathbf{1}_{\{\tau_{d_H}>s\}}h(\widetilde{W}%
_s)e^{\lambda
L^{d_H}_{s}-\alpha_{+}A_{s}^{(d_H,+)}-\alpha_{-}A_{%
s}^{(d_H,-)}}\big]&=&e^{-s\alpha_{-}}\mathbb{E^P}\big[\mathbf{1}%
_{\{M_{s}<d_H\}} h(\widetilde{W}_{s})\big]\indent \text{if}\indent %
d_H>0 \\
&=&e^{-s\alpha_{+}}\mathbb{E^P}\big[\mathbf{1}_{\{I_{s}>d_H\}} h(%
\widetilde{W}_{s})\big]\indent \text{if}\indent d_H<0
\end{eqnarray*}%

\noindent We now need to recall some well-known results on Brownian
motion first-passage times that one may find for instance in
Karatzas and Shreve (1991).

\newtheorem{fpt}[fund]{Lemma}
\begin{fpt}
The following equalities hold for $u>0$ and $a>0$
\begin{equation*}
\mathbb{P}[\tau_a\leq u]=\mathbb{P}[M_u\geq
a]=\frac{2}{\sqrt{2\pi}}\int_{\frac{a}{\sqrt{u}}}^{\infty}e^{-\frac{x^2}{2}}dx
\end{equation*}
Hence, for $u>0$ and $a\in\mathbb{R}$
\begin{equation*}
\mathbb{P}[\tau_a\in du]=\frac{|a|}{\sqrt{2\pi
u^3}}e^{-\frac{a^2}{2u}}du
\end{equation*}
and for $\lambda>0$
\begin{equation*}
\mathbb{E}[e^{-\lambda\tau_a}]=e^{-|a|\sqrt{2\lambda}}
\end{equation*}
where $\tau _{a}=\inf \{t\geq 0;W_{t}=a\}$
\end{fpt}

\newtheorem{mpt}[fund]{Lemma}
\begin{mpt}
For $b\geq0$ and $a\leq b$, the joint density of $(W_u,M_u)$ is
given by :
\begin{equation*}
\mathbb{P}[W_u\in da, M_u\in db]=\frac{2(2b-a)}{\sqrt{2\pi
u^3}}\exp{\big\{-\frac{(2b-a)^2}{2u}\big\}}da~db
\end{equation*}
and likewise, for $b\leq0$ and $a\geq b$ the joint density of
$(W_u,I_u)$ is given by
\begin{equation*}
\mathbb{P}[W_u\in da, I_u\in db]=\frac{2(a-2b)}{\sqrt{2\pi
u^3}}\exp{\big\{-\frac{(2b-a)^2}{2u}\big\}}da db
\end{equation*}
\end{mpt}

\noindent These lemmas provide us with the following interesting
property

\newtheorem{max2}[fund]{Proposition}
\begin{max2}
Let us consider $W_u$ a standard Brownian motion, $I_u$ and $M_u$
respectively its minimum and maximum values up to time $u$.\\ For
any function $h\in L^1(\mathbb{R})$, the quantity
$k_a(u)=\mathbb{E}\big[\mathbf{1}_{\{\tau_a>u\}}h(W_u)\big]$ is
given by
\begin{equation*}
\int_{-\infty}^{\frac{a}{\sqrt{u}}}dv
\frac{e^{-\frac{v^2}{2}}}{\sqrt{2\pi}}h(v\sqrt{u})-\int_{-\infty}^{-\frac{a}{\sqrt{u}}}dv
\frac{e^{-\frac{v^2}{2}}}{\sqrt{2\pi}}h(v\sqrt{u}+2a)\indent
\text{if}\indent a>0
\end{equation*}
\begin{equation*}
\int_{-\infty}^{-\frac{a}{\sqrt{u}}}dv
\frac{e^{-\frac{v^2}{2}}}{\sqrt{2\pi}}h(-v\sqrt{u})-\int_{-\infty}^{\frac{a}{\sqrt{u}}}dv
\frac{e^{-\frac{v^2}{2}}}{\sqrt{2\pi}}h(-v\sqrt{u}+2a)\indent
\text{if}\indent a<0
\end{equation*}
\end{max2}

\begin{proof}
We first observe that
\begin{eqnarray*}
\mathbb{E^{P}}\big[\mathbf{1}_{\{\tau _{a}>u\}}h(W_{u})\big] &=&\mathbb{E^{P}}%
\big[\mathbf{1}_{\{M_{u}<a\}}h(W_{u})\big]\indent\text{if}\indent a>0 \\
&=&\mathbb{E^{P}}\big[\mathbf{1}_{\{I_{u}>a\}}h(W_{u})\big]\indent\text{if}%
\indent a<0
\end{eqnarray*}%
By symmetry, we only need to show the result in the case $a>0$. From
the previous lemma, we can write
\begin{equation*}
\mathbb{E}\big[\mathbf{1}_{\{M_u<a\}}h(W_{u})\big]=\int_{0}^{a}db\int_{-%
\infty }^{b}dxh(x)\frac{2(2b-x)}{\sqrt{2\pi u^{3}}}\exp {\big\{-\frac{%
(2b-x)^{2}}{2u}\big\}}
\end{equation*}%
Finally, we conclude by applying Fubini's theorem.
\end{proof}\newline
\noindent As a consequence, we can now compute the quantity $C_1$

\newtheorem{otmH3}[fund]{Proposition}
\begin{otmH3}
For a call option such that $d_H>0$ or equivalently $H>S_t$, the
quantity $C_1$ is equal to
\[e^{-(r+\alpha_{-})s}G(K,H,S_t,s)\]
where $s=T-t$ and
\begin{eqnarray*}
G(K,H,S_t,s)&=&0\indent
if\indent K\geq H\\
G(K,H,S_t,s)&=&S_t e^{s\frac{(b+\sigma)^2}{2}}N_1-K
e^{s\frac{b^2}{2}}N_2\indent if\indent K<H\\
N_1&=&N(\frac{d_H}{\sqrt{s}}-\sqrt{s}(b+\sigma))-N(\frac{d_K}{\sqrt{s}}-\sqrt{s}(b+\sigma))\\&&-e^{2(b+\sigma)
d_H}\big(N(-\frac{d_H}{\sqrt{s}}-\sqrt{s}(b+\sigma))-N(\frac{d_K-2d_H}{\sqrt{s}}-\sqrt{s}(b+\sigma))\big)\\
N_2&=&N(\frac{d_H}{\sqrt{s}}-\sqrt{s}b)-N(\frac{d_K}{\sqrt{s}}-\sqrt{s}b)\\&&-e^{2b
d_H}\big(N(-\frac{d_H}{\sqrt{s}}-\sqrt{s}b)-N(\frac{d_K-2d_H}{\sqrt{s}}-\sqrt{s}b)\big)
\end{eqnarray*}
where $N(x)=\frac{1}{\sqrt{2\pi}}\int_{-\infty}^{x} dy
e^{-\frac{y^2}{2}}$
\end{otmH3}

\begin{proof}
We apply Proposition II.12 in the case $a>0$ with $h(x)=(S_0
e^{\sigma x}-K)_{+}e^{bx-2\lambda (x-d_H)_{+}}$.
\end{proof}

\newtheorem{itmH3}[fund]{Proposition}
\begin{itmH3}
For a call option such that $d_H<0$ or equivalently $H<S_t$, the
quantity $C_1$ is given by
\[e^{-(r+\alpha_{+})s}J(K,H,S_t,s)\]
where $s=T-t$
\begin{eqnarray*}
J(K,H,S_t,s)&=&S_t
e^{s\frac{(b-2\lambda+\sigma)^2}{2}}N_1(d_1,d_2)-K e^{s\frac{(b-2\lambda)^2}{2}} N_2(d_1,d_2)\\
N_1(d_1,d_2)&=&N(-\frac{d_1}{\sqrt{s}}+\sqrt{s}(b+\sigma-2\lambda))-e^{2(b+\sigma-2\lambda)
d_H}N(\frac{d_2}{\sqrt{s}}+\sqrt{s}(b+\sigma-2\lambda))\\
N_2(d_1,d_2)&=&N(-\frac{d_1}{\sqrt{s}}+\sqrt{s}(b-2\lambda))-e^{2(b-2\lambda)
d_H}N(\frac{d_2}{\sqrt{s}}+\sqrt{s}(b-2\lambda))\\
(d_1,d_2)&=&(d_K,2d_H-d_K)\indent if\indent K>H\\
(d_1,d_2)&=&(d_H,d_H)\indent if\indent K\leq H
\end{eqnarray*}
where $N(x)=\frac{1}{\sqrt{2\pi}}\int_{-\infty}^{x} dy
e^{-\frac{y^2}{2}}$
\end{itmH3}

\begin{proof}
We apply Proposition II.12 in the case $a<0$ with $h(x)=(S_0
e^{\sigma x}-K)_{+}e^{bx-2\lambda (x-d_H)_{+}}$.
\end{proof}

In order to compute $C_2$, it is useful to exhibit a result similar
to the one obtained in
Proposition II.5 to obtain the Laplace transform of the joint density of $%
(B_t,L^a_t,A_t^{(a,+)},A_t^{(a,-)})$.\newline

\newtheorem{excur2}[fund]{Proposition}
\begin{excur2}
Let us consider $W_t$ a standard Brownian motion, $L^a_t$ its local
time at the level $a$, $A_t^{(a,+)}$ and $A_t^{(a,-)}$ respectively
the time spent above and below $a$ by the Brownian motion $W$ until
time $t$.\\ For any function $h\in L^1(\mathbb{R})$, the Laplace
Transform $\int_0^{\infty} dt e^{-\frac{\theta}{2}t}g_a(t)$ of the
quantity
$g_a(t)=\mathbb{E}\big[\mathbf{1}_{\{\tau_a<t\}}h(W_t)\exp(\lambda
L^a_t)\exp(-\mu A_t^{(a,+)}-\nu A_t^{(a,-)})\big]$ is given by
\begin{eqnarray*}
2 e^{-a\sqrt{\theta+2\nu}}\frac{\bigg(\int_0^{\infty}dx
e^{-x\sqrt{\theta+2\mu}}h(a+x)+\int_0^{\infty}dx
e^{-x\sqrt{\theta+2\nu}}h(a-x)\bigg)}{\sqrt{\theta+2\mu}+\sqrt{\theta+2\nu}
-2\lambda}\indent\textrm{if}\indent
a>0\\
2 e^{a\sqrt{\theta+2\mu}}\frac{\bigg(\int_0^{\infty}dx
e^{-x\sqrt{\theta+2\mu}}h(a+x)+\int_0^{\infty}dx
e^{-x\sqrt{\theta+2\nu}}h(a-x)\bigg)}{\sqrt{\theta+2\mu}+\sqrt{\theta+2\nu}
-2\lambda}\indent\textrm{if}\indent a<0
\end{eqnarray*}
for $\theta$ large enough, as seen before.
\end{excur2}

\begin{proof}
Let us prove this result in the case $a>0$; it easily yields to the
case $a<0$.\newline \noindent We first write
\begin{equation*}
g_a(t)=e^{-\nu
t}\mathbb{E}\big[\mathbf{1}_{\{\tau_a<t\}}h(W_t)\exp(\lambda
L^a_t)e^{-(\mu-\nu) A_t^{(a,+)}}\big]
\end{equation*}
Then
\begin{equation*}
I(\theta)=\int_0^{+\infty}dt~e^{-t\frac{\theta+2\nu}{2}}\mathbb{E}\big[%
\mathbf{1}_{\{\tau_a<t\}}h(W_t)\exp(\lambda L^a_t)e^{-(\mu-\nu)
A_t^{(a,+)}}\big]
\end{equation*}

\noindent We now use the strong Markov property and observe that $B_t=W_{t+\tau_a}-W_{%
\tau_a}=W_{t+\tau_a}-a$ is a Brownian motion. Next, we compute the
quantity
\begin{eqnarray*}
\mathbb{E}\big[\mathbf{1}_{\{\tau_a<t\}}h(W_t)\exp(\lambda
L^a_t)e^{-(\mu-\nu)
A_t^{(a,+)}}\big]&=&\mathbb{E}\big[\mathbf{1}_{\{\tau_a<t\}}h(B_{t-\tau_a}+a)%
\exp(\lambda L_{t-\tau_a})e^{-(\mu-\nu) A_{t-\tau_a}^{+}}\big] \\
&=&\int_0^t ds \frac{a e^{-\frac{a^2}{2s}}}{\sqrt{2\pi s^3}} \mathbb{E}\big[%
h(B_{t-s}+a)\exp(\lambda L_{t-s})e^{-(\mu-\nu) A_{t-s}^{+}}\big] \\
&=&\int_0^t ds \frac{a e^{-\frac{a^2}{2(t-s)}}}{\sqrt{2\pi (t-s)^3}} \mathbb{%
E}\big[h(B_s+a)\exp(\lambda L_s) e^{-(\mu-\nu) A_s^{+}}\big]
\end{eqnarray*}
since
\begin{eqnarray*}
L^a_t(a+B_{({\cdot-\tau_a})_{+}})=L_{({t-\tau_a})_{+}} \\
A_t^{(a,+)}=\int_0^t ds \mathbf{1}_{\{B_{(s-\tau_a)_{+}}>0\}}=A^{+}_{({t-\tau_a}%
)_{+}}
\end{eqnarray*}
Hence, applying Fubini's theorem and Proposition II.6 we obtain
\begin{eqnarray*}
I(\theta)&=&\int_0^{\infty} ds e^{-\frac{\theta}{2}s} \mathbb{E}\big[%
h(a+B_s)\exp(\lambda L_s)\exp(-\mu A_s^{+}-\nu
A_s^{-})\big]\int_0^{\infty}
du e^{- \frac{\theta+2\nu}{2}u}\frac{|a| e^{-\frac{a^2}{2u}}}{\sqrt{2\pi u^3}%
} \\
&=&2e^{-a\sqrt{\theta+2\nu}}\frac{\bigg(\int_0^{\infty}dx e^{-x\sqrt{%
\theta+2\mu}}h(a+x)+\int_0^{\infty}dx e^{-x\sqrt{\theta+2\nu}}h(a-x)%
\bigg)}{\sqrt{\theta+2\mu}+\sqrt{\theta+2\nu} -2\lambda}
\end{eqnarray*}
\end{proof}

\newtheorem{otmH}[fund]{Proposition}
\begin{otmH}
In the case $H\leq K$, the Laplace transform with respect to the
variable $T-t$ of the quantity $C_2$ is given by the following
formula:
\[I(\theta)=M(\theta)\frac{N(\theta)}{D(\theta)}\]
where
\[\theta>(\sigma+b-2\lambda)^2-2(r+\alpha_{+})\] and
\begin{eqnarray*}
M(\theta)&=&\bigg(\frac{H}{S_t}\bigg)^{\frac{b-\sqrt{\theta+2(r+\alpha_{-})}}{\sigma}}\indent
if\indent H>S_t\\
M(\theta)&=&\bigg(\frac{S_t}{H}\bigg)^{\frac{2\lambda-b-\sqrt{\theta+2(r+\alpha_{+})}}{\sigma}}\indent
if\indent H<S_t\\
D(\theta)&=&\frac{\sqrt{\theta+2(r+\alpha_{+})}+\sqrt{\theta+2(r+\alpha_{-})}-2\lambda}{2}\\
N(\theta)&=&\frac{H}{\sqrt{\theta+2(r+\alpha_{+})}+2\lambda-\sigma-b}\bigg(\frac{H}{K}\bigg)^\frac{\sqrt{\theta+2(r+\alpha_{+})}+2\lambda-\sigma-b}{\sigma}\\&&-\frac{K}{\sqrt{\theta+2(r+\alpha_{+})}+2\lambda-b}\bigg(\frac{H}{K}\bigg)^\frac{\sqrt{\theta+2(r+\alpha_{+})}+2\lambda-b}{\sigma}
\end{eqnarray*}
\end{otmH}

\begin{proof}
We prove this result by applying Proposition II.8 and Proposition
II.15 and noticing that $(S_0 e^{\sigma (x+d_H)}-K)_{+}=(H e^{\sigma
x}-K)_{+}$
\end{proof}

\newtheorem{itmH}[fund]{Proposition}
\begin{itmH}
In the case $H\geq K$, the Laplace transform with respect to the
variable $T-t$ of the quantity $C_2$ is given by the formula:
\begin{eqnarray*}
I(\theta)&=&M(\theta)\frac{N_1(\theta)+N_2(\theta)}{D(\theta)}\\
M(\theta)&=&\bigg(\frac{H}{S_t}\bigg)^{\frac{b-\sqrt{\theta+2(r+\alpha_{-})}}{\sigma}}\indent
if\indent H>S_t\\
M(\theta)&=&\bigg(\frac{S_t}{H}\bigg)^{\frac{2\lambda-b-\sqrt{\theta+2(r+\alpha_{+})}}{\sigma}}\indent
if\indent H<S_t\\
D(\theta)&=&\frac{\sqrt{\theta+2(r+\alpha_{+})}+\sqrt{\theta+2(r+\alpha_{-})}-2\lambda}{2}\\
N_1(\theta)&=&\frac{H}{\sqrt{\theta+2(r+\alpha_{+})}+2\lambda-\sigma-b}-\frac{K}{\sqrt{\theta+2(r+\alpha_{+})}+2\lambda-b}\\
N_2(\theta)&=&\frac{H}{\sqrt{\theta+2(r+\alpha_{-})}+\sigma+b}\bigg(1-\big(\frac{K}{H}\big)^{\frac{\sqrt{\theta+2(r+\alpha_{-})}+\sigma+b}{\sigma}}\bigg)\\&&
-\frac{K}{\sqrt{\theta+2(r+\alpha_{-})}+b}\bigg(1-\big(\frac{K}{H}\big)^{\frac{\sqrt{\theta+2(r+\alpha_{-})}+b}{\sigma}}\bigg)
\end{eqnarray*}
where $\theta>(\sigma+b-2\lambda)^2-2(r+\alpha_{+})$
\end{itmH}

\begin{proof}
This result is immediately derived from Proposition II.9 and
Proposition II.15.
\end{proof}

\subsection*{E.\indent Extension to a Moving High-Water Mark}

We now wish to take into account the fact that the threshold
triggering the performance fees may accrue at the risk-free rate. As
a consequence, we define $\widetilde{f}$ as
\begin{equation*}
\widetilde{f}(t,S_t)=\mu a \mathbf{1}_{\{S_t>H e^{r t}\}}
\end{equation*}

\newtheorem{fund2}[fund]{Proposition}
\begin{fund2}
There exists a unique solution to the stochastic differential
equation
\begin{equation}
\frac{dS_t}{S_t}=(r+\alpha-c-\widetilde{f}(t,S_t))dt+\sigma dW_t
{\label{fundSDE2}}
\end{equation}
\end{fund2}

\begin{proof}
Let us denote $Y_t=\frac{\ln(S_t e^{-rt})}{\sigma}$. Applying Itô's
formula, we can see that $Y_{t}$ satisfies the following equation
\begin{equation*}
dY_{t}=dW_{t}+\psi (e^{\sigma Y_{t}})dt
\end{equation*}%
where $\psi (x)=-\frac{\sigma ^{2}}{2}+\alpha-c-f(x)$ and $f$
denotes the performance fees function defined in equation
(\ref{perffeesf}).
\newline $\psi $ is Borel locally bounded, consequently we may again
apply Zvonkin theorem that ensures strong existence and pathwise
uniqueness of the solution of ({\ref{fundSDE2}}).
\end{proof}

Let us denote $\widetilde{S}_t=S_te^{-rt}$; we seek to construct a
probability measure $\mathbb{\widehat{Q}}$ under which
\begin{equation*}
\widetilde{S}_t=S_{0}\exp (\sigma \widehat{W}_{t}){\label{rem2}}
\end{equation*}
where $\widehat{W}_{t}$ is a $\mathbb{\widehat{Q}}$ standard
Brownian motion. We briefly extend the results of the previous
section to the case of a moving high-water mark.\newline
\newtheorem{girsanov2}[fund]{Proposition}
\begin{girsanov2}
There exists an equivalent martingale measure $\mathbb{\widehat{Q}}$
under which the Net Asset Value dynamics satisfy the stochastic
differential equation
\begin{equation}
\frac{dS_t}{S_t}=(r+\frac{\sigma^2}{2}) dt +\sigma d\widehat{W}_t
\end{equation}
Moreover,
\begin{equation}
\mathbb{Q}_{|\mathcal{F}_t}=Z_t\cdot\mathbb{\widehat{Q}}_{|\mathcal{F}_t}
\end{equation}
where
\[Z_t=\exp\big(\int_0^t
\big(b-\frac{f(\widetilde{S}_u)}{\sigma}\big)d\widehat{W}_u-\frac{1}{2}\int_0^t
\big(b-\frac{f(\widetilde{S}_u)}{\sigma}\big)^2du\big)\] and
\[b=\frac{\alpha-c-\frac{\sigma^2}{2}}{\sigma}\]
\end{girsanov2}

\newtheorem{expgirs2}[fund]{Lemma}
\begin{expgirs2}
Let us define $d_H$, $\lambda$, $\alpha_{+}$, $\alpha_{-}$ and
$\phi$ as follows:
\[d_H=\frac{\ln(\frac{H}{S_0})}{\sigma},\indent\lambda=\frac{\mu a}{2\sigma}\]
\[\alpha_{+}=2\lambda^2+\frac{b^2}{2}-2\lambda b, \indent\alpha_{-}=\frac{b^2}{2}\]
\[\phi(x)=e^{b x-2\lambda (x-d_H)_{+}}\]
We then obtain:
\begin{equation}
Z'_t=e^{2\lambda(-d_H)_{+}}\phi(\widehat{W}_t)\exp(\lambda
L^{d_H}_t)\exp(-\alpha_{+}A_t^{(d_H,+)}-\alpha_{-}A_t^{(d_H,-)})
\end{equation}
\end{expgirs2}

For the sake of simplicity, we write in this paragraph the strike as
$Ke^{rT}$ and need to compute
\begin{equation}
C(t,S_t)=e^{-r(T-t)}\mathbb{E^Q}\big[(S_T-Ke^{rT})_{+}|\mathcal{F}_t\big]
\end{equation}
The pricing formulas\footnote{All full proofs may be obtained from
the authors.} are derived in a remarkably simple manner by setting
$r=0$ in the results obtained in II.C and II.D.

\section{Numerical Approaches to the NAV option prices}

At this point, we are able to compute option prices thanks to
Laplace Transforms techniques (see Abate and Whitt (1995)) or Fast
Fourier Transforms techniques (see Walker (1996)) . We can observe
that if Monte Carlo simulations were performed in order to obtain
the NAV option price, the number of such simulations would be fairly
large because of the presence of an indicator variable in the Net
Asset Value dynamics. The computing time involved in the inversion
of Laplace transforms is remarkably lower compared to the one
attached to Monte Carlo simulations. The times to maturity
considered below are chosen to be less or equal to one year in order
to avoid the high water mark reset arising for more distant
maturities. Taking into account the reset feature would lead to
computations analogous to the ones involved in forward start options
and is not the primary focus of this paper.\newline Tables 1 to 4
show that the call price is an increasing function of the excess
performance $\alpha$, and in turn drift $\mu$, a result to be
expected.\newline
\newline The call price is also increasing with the high water mark level $H$ as
incentive fees get triggered less often.
\newline
\newline Table 5 was just meant
to check the exactitude of our coding program : by choosing $a=0$
and $\alpha=0$, the NAV call option pricing problem is reduced to
the Merton (1973) formula. Table 5 shows that the prices obtained by
inversion of the Laplace transform are remarkably close to those
provided by the Merton analytical formula.

\begin{equation*}
\begin{tabular}{|c|cc|}
\multicolumn{2}{c}{Table 1}\\
\multicolumn{2}{c}{Call Option Prices at a volatility level $\sigma=20\%$} \\
\multicolumn{2}{c}{$H=\$ 85$, $S_0=\$ 100$, $\alpha=10\%$, $r=2\%$,
$c=2\%$, $a=20\%$, $\mu=15\%$} \\\hline Strike / Maturity & 6 Months
& 1 Year \\\hline $90\%$ & $\$14.5740$ & $\$18.9619$
\\\hline $100\%$ & $\$7.6175$ & $\$12.1470$ \\\hline $110\%$ &
$\$3.3054$ &$\$7.2058$
\\\hline
\multicolumn{2}{c}{$H=S_0=\$ 100$, $\alpha=10\%$, $r=2\%$, $c=2\%$,
$a=20\%$ and $\mu=15\%$} \\\hline Strike / Maturity & 6 Months & 1
Year \\\hline $90\%$ & $\$15.0209$ & $\$19.6866$ \\\hline $100\%$ &
$\$7.8346$ & $\$12.5922$ \\\hline $110\%$ & $\$3.3837$ & $\$7.4427$
\\\hline
\multicolumn{2}{c}{$H=\$ 115$, $S_0=\$ 100$, $\alpha=10\%$, $r=2\%$,
$c=2\%$, $a=20\%$, $\mu=15\%$} \\\hline Strike / Maturity & 6 Months
& 1 Year \\\hline $90\%$ & $\$15.7095$ & $\$20.8464$ \\\hline
$100\%$ & $\$8.4147$ & $\$13.5815$ \\\hline $110\%$ & $\$3.7084$ &
$\$8.1198$ \\\hline
\end{tabular}%
\end{equation*}
\begin{equation*}
\begin{tabular}{|c|cc|}
\multicolumn{2}{c}{Table 2} \\
\multicolumn{2}{c}{Call Option Prices at a volatility level $\sigma=20\%$} \\
\multicolumn{2}{c}{$H=\$ 85$, $S_0=\$ 100$, $\alpha=15\%$, $r=2\%$,
$c=2\%$, $a=20\%$, $\mu=20\%$} \\\hline Strike / Maturity & 6 Months
& 1 Year \\\hline $90\%$ & $\$16.3804$ & $\$22.6562$ \\\hline
$100\%$ & $\$8.9668$ & $\$15.1925$ \\\hline $110\%$ & $\$4.1091$ &
$\$9.4795$ \\\hline \multicolumn{2}{c}{$H=S_0=\$ 100$,
$\alpha=15\%$, $r=2\%$, $c=2\%$, $a=20\%$ and $\mu=20\%$} \\\hline
Strike / Maturity & 6 Months & 1 Year \\\hline $90\%$ & $\$16.9611$
& $\$23.6036$ \\\hline $100\%$ & $\$9.2703$ & $\$15.8190$ \\\hline
$110\%$ & $\$4.2276$ & $\$9.8398$ \\\hline \multicolumn{2}{c}{$H=\$
115$, $S_0=\$ 100$, $\alpha=15\%$, $r=2\%$, $c=2\%$, $a=20\%$,
$\mu=20\%$} \\\hline Strike / Maturity & 6 Months & 1 Year \\\hline
$90\%$ & $\$17.9362$ & $\$25.2503$ \\\hline $100\%$ & $\$10.1156$ &
$\$17.2719$ \\\hline $110\%$ & $\$4.7300$ & $\$10.8943$ \\\hline
\end{tabular}%
\end{equation*}

\begin{equation*}
\begin{tabular}{|c|cc|}
\multicolumn{2}{c}{Table 3} \\
\multicolumn{2}{c}{Call Option Prices at a volatility level $\sigma=40\%$} \\
\multicolumn{2}{c}{$H=\$ 85$, $S_0=\$ 100$, $\alpha=10\%$, $r=2\%$,
$c=2\%$, $a=20\%$, $\mu=15\%$} \\\hline Strike / Maturity & 6 Months
& 1 Year \\\hline $90\%$ & $\$18.8245$ & $\$25.3576$ \\\hline
$100\%$ & $\$13.2042$ & $\$19.9957$ \\\hline $110\%$ & $\$8.9804$ &
$\$15.6276$ \\\hline \multicolumn{2}{c}{$H=S_0=\$ 100$,
$\alpha=10\%$, $r=2\%$, $c=2\%$, $a=20\%$ and $\mu=15\%$} \\\hline
Strike / Maturity & 6 Months & 1 Year \\\hline $90\%$ & $\$19.1239$
& $\$25.8231$ \\\hline $100\%$ & $\$13.3979$ & $\$20.3534$ \\\hline
$110\%$ & $\$9.1012$ & $\$15.8949$ \\\hline \multicolumn{2}{c}{$H=\$
115$, $S_0=\$ 100$, $\alpha=10\%$, $r=2\%$, $c=2\%$, $a=20\%$,
$\mu=15\%$} \\\hline Strike / Maturity & 6 Months & 1 Year \\\hline
$90\%$ & $\$19.5128$ & $\$26.4273$ \\\hline $100\%$ & $\$13.7277$ &
$\$20.8726$ \\\hline $110\%$ & $\$9.3409$ & $\$16.3134$ \\\hline
\end{tabular}%
\end{equation*}

\begin{equation*}
\begin{tabular}{|c|cc|}
\multicolumn{2}{c}{Table 4}\\
\multicolumn{2}{c}{Call Option Prices at a volatility level $\sigma=40\%$} \\
\multicolumn{2}{c}{$H=\$ 85$, $S_0=\$ 100$, $\alpha=15\%$, $r=2\%$,
$c=2\%$, $a=20\%$, $\mu=20\%$} \\\hline Strike / Maturity & 6 Months
& 1 Year \\\hline $90\%$ & $\$20.3926$ & $\$28.6499$ \\\hline
$100\%$ & $\$14.4928$ & $\$22.8861$ \\\hline $110\%$ & $\$9.9903$ &
$\$18.1179$ \\\hline \multicolumn{2}{c}{$H=S_0=\$ 100$,
$\alpha=10\%$, $r=2\%$, $c=2\%$, $a=20\%$ and $\mu=15\%$} \\\hline
Strike / Maturity & 6 Months & 1 Year \\\hline $90\%$ & $\$20.7978$
& $\$29.2995$ \\\hline $100\%$ & $\$14.7618$ & $\$23.3938$ \\\hline
$110\%$ & $\$10.1615$ & $\$18.5044$ \\\hline
\multicolumn{2}{c}{$H=\$ 115$, $S_0=\$ 100$, $\alpha=10\%$, $r=2\%$,
$c=2\%$, $a=20\%$, $\mu=15\%$} \\\hline Strike / Maturity & 6 Months
& 1 Year \\\hline $90\%$ & $\$21.3417$ & $\$30.1555$ \\\hline
$100\%$ & $\$15.2260$ & $\$24.1402$ \\\hline $110\%$ & $\$10.5042$ &
$\$19.1158$ \\\hline
\end{tabular}%
\end{equation*}

\begin{equation*}
\begin{tabular}{|c|cc|cc|}
\multicolumn{5}{c}{Table 5}\\
\multicolumn{5}{c}{NAV Call Option Prices when $\mu=0$ at a volatility level $\sigma=40\%$} \\
\multicolumn{5}{c}{$S_0=\$ 100$, $r=2\%$, $c=0.3\%$}\\
\hline Maturity & \multicolumn{2}{|c|}{6 months} & \multicolumn{2}{|c|}{1 year} \\
Strike & Laplace Transform & Merton formula & Laplace Transform &
Merton formula \\\hline \multicolumn{1}{|c|}{90\%} &
\multicolumn{1}{|c|}{\$12.3324} & \multicolumn{1}{|c|}{\$12.3324} &
\multicolumn{1}{|c|}{\$14.577} & \multicolumn{1}{|c|}{\$14.577}
\\\hline \multicolumn{1}{|c|}{100\%} &
\multicolumn{1}{|c|}{\$6.0375} & \multicolumn{1}{|c|}{\$6.0375} &
\multicolumn{1}{|c|}{\$8.7434} & \multicolumn{1}{|c|}{\$8.7434}
\\\hline \multicolumn{1}{|c|}{110\%} &
\multicolumn{1}{|c|}{\$2.4287} & \multicolumn{1}{|c|}{\$2.4287} &
\multicolumn{1}{|c|}{\$4.8276} &
\multicolumn{1}{|c|}{\$4.8276}\\\hline
\end{tabular}
\end{equation*}

\pagebreak
\section{Conclusion}

In this paper, we proposed a pricing formula for options on hedge
funds that accounts for the high-water mark rule defining the
performance fees paid to the fund managers. The geometric Brownian
motion dynamics chosen for the hedge fund Net Asset Value allowed us
to exhibit an explicit expression of the Laplace transform in
maturity of the option price through the use of Brownian local
times. Numerical results obtained by inversion of these Laplace
transforms display the influence of key parameters such as
volatility or moneyness on the NAV call price.

\pagebreak

\pagebreak

\section*{Appendix : Excursion Theory}

\begin{proof}
\textbf{of Proposition II.5}: We use the Master formula exhibited in
Brownian excursion theory; for more details see Chapter XII in Revuz
and Yor (2005) the notation of which we borrow:\newline
\noindent $n$ denotes the It\^{o} characteristic measure of excursions and $%
n_{+} $ is the restriction of $n$ to positive excursions;\newline
\noindent $V(\epsilon)=\inf\{t>0; \epsilon(t)=0\}$ for
$\epsilon\in\mathbf{W}_{exc}$ the space of excursions,\newline
\noindent $(\tau_l)_{l\geq 0}$ is the inverse local time of the
Brownian motion.\newline
\newline

\noindent We can write
\begin{equation*}
\mathbb{E}\bigg[\int_0^{\infty} dt e^{-\frac{\theta}{2}t}
h(W_t)\exp(\lambda L_t)\exp(-\mu A_t^{+}-\nu A_t^{-})\bigg]=I\cdot J
\end{equation*}
where
\begin{eqnarray*}
I&=&\mathbb{E}\big[\int_0^{\infty} dl e^{-\frac{\theta}{2}\tau_l}
e^{\lambda l} \exp(-\mu A_{\tau_l}^{+}-\nu A_{\tau_l}^{-}) \big] \\
&=&\int_0^{\infty}dl \exp\bigg(l\big(\lambda -\int n(d\epsilon) (1-e^{-\frac{%
\theta}{2}V-\mu A_V^{+}-\nu A_V^{-}})\big)\bigg) \\
&=&\frac{1}{\int n(d\epsilon) \big(1-e^{-\frac{\theta}{2}V-\mu
A_V^{+}-\nu
A_V^{-}}\big)-\lambda} \\
&=&\frac{1}{\frac{\sqrt{\theta+2\mu}+\sqrt{\theta+2\nu}}{2}-\lambda}
\end{eqnarray*}

\noindent and
\begin{equation*}
J=\int_0^{\infty} ds e^{-\frac{\theta}{2}s}\bigg\{e^{-\mu s} n_{+}\big(%
h(\epsilon_s)\mathbf{1}_{\{s<V\}}\big)+e^{-\nu s} n_{+}\big(h(-\epsilon_s)%
\mathbf{1}_{\{s<V\}}\big) \bigg\}
\end{equation*}

\noindent Next, we use the result
\begin{equation}
n_{+}\big(\epsilon_s\in dy;s<V\big)=\frac{y}{\sqrt{2\pi s^3}}e^{-\frac{y^2}{2s}}dy\indent(y>0) {%
\label{den}}
\end{equation}
\noindent and obtain
\begin{equation}
J=\int_0^{\infty}dx e^{-x\sqrt{\theta+2\mu}}h(x)+\int_0^{\infty}dx
e^{-x\sqrt{\theta+2\nu}}h(-x){\label{res3}}
\end{equation}
\noindent where the proof of equation (\ref{res3}) comes from the fact that in (\ref%
{den}) the density of $n_{+}$ as a function of $s$, is precisely the density of $%
T_y=\inf\{t: B_t=y\}$, and $\mathbb{E}[e^{-\lambda T_y}]=e^{-y\sqrt{2\lambda}%
}$.
\end{proof}\\
This example of application of excursion theory is one of the
simplest illustrations of Feynman-Kac type computations which may be
obtained with excursion theory arguments. For a more complete story,
see Jeanblanc, Pitman and Yor (1997).

\end{document}